\documentclass[12pt]{article}
\usepackage{graphics,latexsym,amssymb,amsmath,amscd}
\RequirePackage{graphics}
\usepackage{multirow}
\usepackage[all,cmtip]{xy}
\usepackage{xcolor}
\def\[#1\]{\begin{eqnarray*}#1\end{eqnarray*}}

\def\phi{\varphi}

\input xy
\xyoption{all}
\newtheorem{thm}{Theorem}[section]
\newtheorem{dfn}[thm]{Definition}
\newtheorem{rmk}[thm]{Remark}
\newtheorem{cor}[thm]{Corollary}
\newtheorem{prop}[thm]{Proposition}
\newtheorem{lemma}[thm]{Lemma}

\newcommand{\Pf}{{\em Proof}. }
\newcommand{\EPf}{\hbox{}\hfill$\Box$\vspace{.5cm}}

\newcommand{\R}{{{\mathbb R}}}

\def\phi{\varphi}

\def\N{{\mathbb N}}

\def\R{{\mathbb R}}

\date{}
\title{Self-adjoint extensions of singular Sturm-Liouville operators on graphs and Weyl's law}
\author{Elisha Falbel}
 
\begin{document}

\maketitle
\newcommand{\D}{\mbox{$\cal D$}}
\newtheorem{df}{Definition}[section]
\newtheorem{te}{Theorem}[section]
\newtheorem{co}{Corollary}[section]
\newtheorem{po}{Proposition}[section]
\newtheorem{lem}{Lemma}[section]
\newcommand{\Ad}{\mbox{Ad}}
\newcommand{\ad}{\mbox{ad}}
\newcommand{\im}[1]{\mbox{\rm im\,$#1$}}
\newcommand{\bm}[1]{\mbox{\boldmath $#1$}}
\newcommand{\sime}{\mbox{sim}}

\begin{abstract}
We study self-adjoint extensions of a second order differential operator of Sturm-Liouville type on a graph.    We relate self-adjointness of the operator to the existence of non-complete trajectories  of the Hamiltonian vector field defined by its principal symbol outside the vertices.  We define Kirchhoff conditions at the vertices which guarantee a self-adjoint extension analogous to the case of quantum graphs.  The singular vertices may be interpreted as introducing a singular potential at those points.  We also establish a Weyl's law for the spectrum asymptotics. 

\end{abstract}

\section{Introduction}


In \cite{T}(2020,2024) and \cite{dVLB}(2022) it is proven that the self-adjointness of certain classes of symmetric ordinary differential operators (more generally, pseudo-differential operators in \cite{T}) on the circle occurs if and only if the Hamiltonian vector field defined on the cotangent bundle by the principal symbol is complete.  They state it as the equivalence between quantum completeness (essential self-adjointness of an operator) and classical completeness (completeness of the Hamiltonian vector field associated to the principal symbol of an operator).  The simplest situation is that of real symmetric second-order differential operators, that is Sturm-Liouville operators, with smooth coefficients on the circle.
 
 A basic observation is that even for the singular Sturm-Liouville operator defined on an interval, the existence of non-complete trajectories at the boundary of the interval counts the number of self-adjoint extensions.  This follows from 
 Proposition \ref{proposition:Hamiltonian} and the criterium in Theorem \ref{theorem:deg-LCP}.
 
 Instead of working on a circle, we might patch together several intervals and their corresponding singular Sturm-Liouville operators to 
 compute that the self-adjoint extensions of a second order real differential operator  on a graph.

 We are first interested  to understand how incompleteness is related to self-adjoint extensions of Sturm-Liouville operators.  
 
 \vspace{.5cm}

{\bf Theorem A}
Consider a symmetric differential operator defined on a finite graph $\mathcal{G}$ given by the expression
$$
Py=-(py')'+qy, 
$$
where $p$ and $q$ are continuous functions (smooth up to the vertices)  and $p$  has non-degenerate zeros and is  non-vanishing on the complement of the vertices $\mathcal{V}$.  Suppose that there are $N$ endpoints of edges where $p$ is either non-vanishing or has a simple zero.   Then,
\begin{enumerate}
\item
The self-adjoint extensions of  the operator $P$ on $L^2(\mathcal{G})$ defined on smooth functions with compact support on the complement of the vertices are parametrized by $U(N)$. 
\item 
The closure of $P$ is essentially self-adjoint if and only if  the Hamiltonian vector field corresponding to the Hamiltonian defined
by the the principal symbol $-p(x)\xi^2$ on the cotangent bundle of $\mathcal{G} \setminus \mathcal{V}$ is complete. 
\end{enumerate}

  In particular, if at all endpoints  adjacent to the same  vertex $p$ has the same order,  the self-adjoint extensions of a second order real symmetric operator on the circle are parametrized by the unitary group $U(2N)$ where $N$ is the number of vertices where $p$ vanishes with order one.

  The multiple interval Sturm-Liouville problem where the operator is defined on functions of compact support on the complement of the singular set was analyzed before in \cite{EZ} and \cite{KZ}  in a more general case where the coefficients are not necessarily smooth.  They describe explicit self-adjoint extensions using boundary conditions.
We also observe that the number of singular points in the theorem could be infinite.  This follows from a generalization of Sturm-Liouville equations with an infinite number of singularities which goes back to 
\cite{BF} and \cite{GK} (see also \cite{EZ1}).  Singular quantum operators  appeared also in the study of contact geometry (see \cite{ABPT}).


  In the special case of the circle one can state the following corollary.
 \vspace{.5cm}
 
{\bf Corollary}
Consider a symmetric differential operator defined on $S^1$ given by the expression
$$
Py=-(py')'+qy, 
$$
where $p$ and $q$ are smooth functions and $p$ has a finite number of zeros which are all non-degenerate and such that $n$ of them are simple.  Then,
\begin{enumerate}
\item
The self-adjoint extensions of the closure of the operator on $L^2(S^1)$ defined on smooth functions with compact support on the complement of the zeros are parametrized by $U(2n)$. 
\item
The self-adjoint extensions of the closure of the operator on $L^2(S^1)$ defined on smooth functions with compact support are parametrized by $U(n)$.
\item 
The number of vertical fibers containing non-complete trajectories of the Hamiltonian vector field corresponding to the Hamiltonian defined
by the the principal symbol $-p(x)\xi^2$ on the cotangent bundle is equal to $n$.

\end{enumerate}

\vspace{.5cm}

The proof is obtained directly from Corollary \ref{corollary:circle} and Proposition \ref{proposition:Hamiltonian}. It is important to note that the vertical fibers are approached 
by  non-complete trajectories when they are above a simple singularity even if the symbol is not smooth at the fibers and the 
Hamiltonian vector field is not defined on the fibers. Indeed, Proposition \ref{proposition:Hamiltonian} describes the situation in an open interval and analyzes the Hamiltonian vector field near the boundary.

In section \ref{section:SLoperators} we define singular quantum graphs (or Sturm-Liouville graphs)  and 
describe a subset of self-adjoint extensions (see section \ref{theorem:singular-quantum-graph}) using local boundary conditions which are the analogs of Kirchhoff conditions in the case of regular quantum graphs.

%
%
%
Given a Sturm-Liouville operator $P$ on an interval $(a,b)$ one defines an endpoint to be of limit circle type if all solutions to the equation $Py=0$ are in $L^2$ near the endpoint (we say moreover that the point is regular if the solutions are bounded).  Otherwise, the point is said to be of limit point type.  
We study then
Sturm-Liouville graphs which do not have LP vertices and such that the principal symbol is positive in each interval.  This is a natural condition so that the spectrum of all self-adjoint extensions be purely discrete.  We obtain Weyl's law for such graphs (see Theorem \ref{theorem:Weyl}) generalizing the formula of \cite{BK} for regular quantum graphs.
In the following, for a selfadjoint operator $P$, we denote
$$
 N(\lambda,P)=\#\{ \ \mbox{eigenvalues of}\ P\leq \lambda\ \}.
$$

\vspace{.5cm}

{\bf Theorem B}

Let $(\mathcal{G},P_{min})$ be Sturm-Liouville graph with vertices which are either regular or of LC type. We suppose that the principal symbol is positive in each interval. Then for any self adjoint extension $P$ of $P_{min}$, we have
$$
 N(\lambda,P)\sim \left(\frac{1}{\sqrt{\pi}}\left( \sum_{e\in R}{L_e}^{1/2}
 +{\sqrt{2}}\sum_{e\in RLC}{L_e}^{1/4}\right) +
  \#\{\ e\in E\ |\ e\in LC\ \}
  \right)\sqrt{\lambda},
$$
where the first sum is made over the edges with both endpoints regular, the second with one regular  endpoint and the other LC and the last term is the number of edges with both endpoints LC (non-regular).

\vspace{.5cm}

I thank F. Naud for introducing me  to the references \cite{dVLB} and \cite{T} and  P.  Dingoyan and F. Naud for innumerable fruitful discussions.

\section{Preliminaries}

For a thorough introduction to self-adjoint extensions of differential operators see \cite{N}, \cite{RS} or \cite{W}.

\subsection{Self-adjoint extensions of symmetric operators}

Here we recall von Neumann extension theory which parametrizes self-adjoint extensions of a symmetric operator by a unitary group. One can consult  \cite{N} (section 14) or \cite{RS} (Chapter X). 

Let $A$ be a closed symmetric operator and $A^*$ its adjoint.  Define the deficiency subspaces
$$
N_+= \ker (A^*-iI)
$$
and
$$
N_-= \ker (A^*+iI).
$$
The dimensions $d_+=\dim N_+$ and $d_-=\dim N_-$ are called deficiency indices of $A$
We have that
$$
D(A^*)=D(A)\oplus N_+\oplus N_-.
$$
The operator $A$ has self-adjoint extensions if and only if $d_+$=$d_-$.  Moreover, $d_=d_-=0$ if and only if  the operator is self-adjoint.

\begin{thm} A self-adjoint extension $S$ of a closed symmetric operator $A$ is given by 
an isometric operator $U$ from  $N_+$ into $N_-$.
Its domain is 
$$
D_{S}=\{ x+ z-Uz\ |\ x\in D_{A}, z\in N_+\ \}
$$
and $P(x+ z-Uz)= A(x)+ i z-i U(z)$.
\end{thm}

The set of self-adjoint extensions is parametrized by the unitary group $U(d)$,
where $d$ is the deficiency index of $A$.

\section{Sturm-Liouville operators on an interval}

See \cite{N}, \cite{Z} or \cite{W}.

\subsection{The minimal and maximal operator, endpoints in the LC or LP case}

  For simplicity sake we will suppose that the 
the operator has continuous symbol and later we will assume that the symbol is smooth.

Consider the differential operator defined on $I=(a,b)$ by the expression
$$
Py=-(py')'+qy, 
$$
where $p$ and $q$ are real continuous functions up to the boundary with $p>0$ on $I$.

The natural Hilbert space is $L^2(I)$ and one defines a maximal operator $P_{max}$ (which is closed) whose domain is
$$
D_{max}=\{ \ f \in L^2(I)\ |  \ Pf\in  L^2(I)\ \}.
$$
For $f\in D_{max}$, we have $f, f'\in AC_{loc}(I)$.

The minimal operator $P_{min}$ is defined to be the closure of the  operator $P$ defined on functions in $D_{max}$ of compact support.  The basic relation between $P_{max}$ and $P_{min}$ is that
$$
P_{max}^*=P_{min}\ \mbox{and}\ P_{min}^*=P_{max}.
$$

\begin{dfn} The deficiency indices of $P$ are defined to be 
$$
d^{\pm}=\dim \ker {(P_{max}\mp i Id)}.
$$
\end{dfn}

We always have $d^-=d^+$ and the common value $d$ satisfies $0\leq d\leq 2$.  Observe also that
$$
\dim D_{max}/D_{min}=2d.
$$

The following classification of endpoints dates back to Weyl (1910).

\begin{dfn}
\begin{enumerate}
\item The endpoint $c=a$ is regular if $1/p\in L^2(a,c)$ for some $c\in (a,b)$.  
\item The endpoint $c=a$  is a LC case (limit circle) if all solutions of the equation $Pu=0$ are in $L^2(a,c)$ for some $c\in (a,b)$.
 \item The endpoint $c=a$  is a LP case (limit point) if it is not LC case.
\end{enumerate}
\end{dfn}

An analogous definition is obtained for the classification of the endpoint $b$.
\vspace{.5cm}

We will use repeatedly the following theorem. It gives a practical criterium to decide if an endpoint is a LP or LC case.  See the appendix for a proof.

\begin{thm} \label{theorem:deg-LCP}
Let $p(x)=(b-x)^n\phi(x)$ near the the endpoint $b$ with $n\in \N$ and $\phi(b)\neq 0$ a $C^1$ function up to $b$.  Then the  endpoint $b$   is a limit point case if and only if $n\geq 2$.
\end{thm}


\subsection{Completeness of the Hamiltonian vector field and the limit point case}\label{section:completeness}

We assume in this section that  the differential operator $L$ is defined on $I=(a,b)$ by the expression
$$
Py=-(py')'+qy, 
$$
where $p$ and $q$ are real smooth functions up to the boundary with $p>0$ on $I$.

The principal symbol of $P$, defined on $T^*I=\{\ (x,\xi dx)\ | \ x\in I, \ \xi\in \R\ \}\simeq I\times \R$ is
$$
s(x,\xi)= -p(x)\xi^2.
$$
The canonical symplectic form on $T^*I$ is given by $\omega=dx\wedge d\xi$ and the Hamiltonian vector field 
associated to the principal symbol (that is, satisfying $\iota_X \omega=ds$) is 
$$
X=\frac{\partial s}{\partial \xi} \frac{\partial }{\partial x}-\frac{\partial s}{\partial x} \frac{\partial }{\partial \xi}=
-2p(x)\xi \frac{\partial }{\partial x}+p'(x)\xi^2 \frac{\partial }{\partial \xi}.
$$

The Hamiltonian flow preserves the level sets $L_E$ of the principal symbol which are described  fixing  $E\in \R^+$ defined as
$$L_E=\{\ (x,\xi)\ | \
-p(x)\xi^2=-E\ \}.
$$
Therefore the flow defined by the solutions to the system
$$
\frac{dx}{dt}=-2p(x)\xi, \ \  \frac{d\xi}{dt}=p'(x)\xi^2 
$$
is obtained solving the differential equation
$$
\frac{dx}{dt}=\mp 2\sqrt{E}\sqrt{p(x)}.
$$
Observe that there are two equations according to the sign of $\xi$.

\begin{dfn}
We will say that the Hamiltonian vector field is complete at the endpoint $b$ (respectively $a$) if all solutions $x(t)$ to the equation above 
satisfy  $|t|\to \infty$ when $x\to b$ (respectively $x\to a$).
\end{dfn}

Clearly, the Hamiltonian vector field is complete if and only if it is complete at both
endpoints.

\begin{prop} \label{proposition:Hamiltonian}
The Hamiltonian vector field $X$ is complete at the endpoint $b$  (respectively at $a$) if and only if for $c\in I$, $\frac{1}{\sqrt{p}}\notin L^1(c,b)$ ( respectively $\frac{1}{\sqrt{p}}\notin L^1(a,c)$).
\end{prop}

\Pf
From the equation $
\frac{dx}{dt}=\mp 2\sqrt{E}\sqrt{p(x)}
$ we obtain
$$
\int_{c}^{x(t)}\frac{1}{\sqrt{p(x)}}dx=\mp 2\sqrt{E}\int_{t_0}^t dt.
$$
Therefore, time  tends to infinity when approaching $b$ if and only if $\frac{1}{\sqrt{p}}\notin L^1(c,b)$ and similarly for the endpoint $a$.
\EPf

\section{Sturm-Liouville on direct sums}
Consider now the differential operator defined on a disjoint union of open intervals
$$
\mathcal{I}=\bigcup_{k=1}^{N}{I_k}.
$$ given by the expression
$$
Py=-(py')'+qy, 
$$
where $p$ and $q$ are real continuous functions up to the boundary on each $I_k$. 

The operator $P$ is defined as a Sturm-Liouville operator on $L^2(I)=\bigoplus_{k=1}^N L^2(I_k)$ with restrictions to each 
$I_k$ denoted by $P_k$.

\vspace{1cm}

\subsection{Maximal and minimal operators on direct sums}

  A natural minimal operator is the closure of $P$ acting on smooth functions with compact support contained in $\bigcup_{k=1}^{N}{I_k}$.
We call $P_{Min}$ this operator.   We have
$$
P_{Min}= \sum_{k=1}^N{P_k}_{_{min}}\circ \pi_k,
$$
where ${P_k}_{_{min}}$ is the minimal operator on the interval $I_k$ as in the previous section and $\pi_k$ is the projection of $L^2(I)$ onto $L^2(I_k)$, for $1\leq k\leq N$.  The operator $P_{Min}$ ignores the relative positions of the intervals $I_k$. is best described as 
$$
\oplus_{k=1}^N {P_k}_{_{min}}.
$$

The natural maximal operator $P_{Max}$ has domain 
$$
D_{Max}=\{ \ f \in L^2(I)\ |\  \mbox{for all}\ k,\ Pf_{\vert_{I_k}}\in  L^2(I_k)\ \}.
$$
That is,
$$
D_{Max}= \bigoplus_{k=1}^{N} D_{max}(I_k),
$$
where $D_{max}(I_k)$ is the maximal operator defined in the interval case.


\vspace{1cm}

The fundamental properties of the maximal and minimal operators are easily obtained and 
are stated in the following lemma (see \cite{EZ} or \cite{Z}, lemma 13.3.1):

\begin{lemma} For the operators $P_{Min}$ and $P_{Max}$ defined as above we have
\begin{enumerate}
\item $P^*_{Min}=P_{Max}$ and $P^*_{Max}=P_{Min}$
\item $P^*_{Min}$ is  densely defined and its deficiency indices are
$$
def^+=def^-=\sum_{k=1}^N def_k,
$$
where $def_k$ is the deficiency index of each ${P_k}_{_{min}}$.
\end{enumerate}
\end{lemma}

\begin{thm}
Consider a symmetric differential operator defined by $\oplus_{k=1}^N {P_k}_{_{min}}$ as above by an expression
$$
Py=-(py')'+qy, 
$$
where $p$ and $q$ are continuous functions (put conditions) and such that $p$ which is non-vanishing on each open interval $I_k$.  Suppose that there are $L$ endpoints of edges $I_k$ which are LC (or regular).   Then,
\begin{enumerate}
\item
The self-adjoint extensions of the closure of the operator on $L^2(\bigcup_{k=1}^{N}{I_k})$ defined on smooth functions with compact support on $\bigcup_{k=1}^{N}{I_k}$ is parametrized by $U(2L)$. 
\item If we suppose, moreover,  that $p$ is smooth, 
then the closure of $P$ is essentially self-adjoint if and only if  the Hamiltonian vector field corresponding to the Hamiltonian defined
by the the principal symbol $-p(x)\xi^2$ on the cotangent bundle of $\bigcup_{k=1}^{N}{I_k}$ is complete. 

\item In particular, suppose $p$ and $q$ are continuous functions and such that $p$ is smooth up to the vertices, with non-degenerate zeros, which is non-vanishing on the complement of the endpoints.  Then $L$ above is the number of endpoints  where $p$ is either non-vanishing or has a simple zero. 
\end{enumerate}

\end{thm}

\vspace{1cm}

\begin{df}[Lagrangian form]\label{definition:lagrange-bilinear}
Given two functions $y,z\in D_{Max}$ define the skew-hermitian bilinear form
$$
[y,z]=yp\bar{z'}-\bar{z}py'.
$$
\end{df}
Observe that, outside the endpoints, $[y,z]$ is absolutely continuous 
and one verifies the following identity:
$$
\bar{z}P y-y\overline{Pz}= [y,z]'.
$$

We have therefore:
\begin{lem} (Green's formula)
For $a_k<\alpha<\beta<b_k$ ($1\leq k\leq N$) and $y,z\in D_{Max}$,
$$
\int_\alpha^\beta(\bar{z}P y-y\overline{Pz})=[y,z](\beta)-[y,z](\alpha).
$$
\end{lem}
Green's formula implies that, for any $y,z\in D_{Max}$, $[y,z]$ has an extension to the end points of $I_k$ and 
$$
\int_{a_k}^{b_k}(\bar{z}P y-y\overline{Pz})=[y,z](b_k)-[y,z](a_k).
$$


\begin{lem}[see  \cite{N,Z}]
If $a_k$ and $b_k$ are in the LP case, then $[y,z](b_k)=[y,z](a_k)=0$ for all $y,z\in D_{Max}$.
\end{lem}

Taking $y\equiv 1$, the constant function, we obtain that, for all $z\in D_{Max}$, 
$pz'(a_k)= pz'(b_k)=0$ if $a_k$ and $b_k$ are in the LP case.

\subsection{Sturm Liouville operators with interior singular points : More constrained maximal and minimal operators}

The maximal and minimal operators defined so far don't depend on the relative position of the intervals on the real line.  Indeed, they are defined as a direct sum of operators. We define now operators which are sensitive to the precise position of the intervals.  Let again
$$
\mathcal{I} =\bigcup_{k=1}^{N}{I_k}.
$$
But now, suppose that the endpoints are identified so that $I=\bigcup_{k=1}^{N}\overline{I_k}$ (with identifications) is a one dimensional manifold.
Define the  operator $P_{max}$ with domain
$$
D_{max}=\{ \ f \in L^2(I)\ |  \ Pf\in  L^2(I)\ \}
$$
and such that $P_{max}$ is  $P_{Max}$ restricted to $D_{max}$.

{

 \begin{rmk} 
Given a distribution $u\in D'(I)$ such that the distributional derivative satisfies $\partial u\in L^1_{loc}(I)$ then $u\in L^1_{loc}(I)$ (it is actually locally absolutely continuous).  We conclude that $pf'$ is locally absolutely continuous and $f$ is locally absolutely continuous where $p$ is nonzero.  
\end{rmk}
}
\vspace{1cm}

Observe that 
$$
D_{max}\subset  D_{Max}\subset L^2(I).
$$

Define $P_{min}$ the closure of the operator defined on the domain 
$$
D'_{min}= \{ \ f \in L^2(I)\ |\ f\in C^{\infty}_0(I)\ \},
$$
with $P_{min}(f)=Pf$ and denote the domain of the closure  $ D_{min}$.

The following lemma is a generalization of the one interval case (see \cite{GK}, \cite{EZ1} and also lemma 2.2 in \cite{dVLB}).

\begin{lemma} Given $P_{min}$ defined as above, the domain of the adjoint $P^*_{min}$ is 
$$
D^*_{min} =D_{max}=\{ \ g \in L^2(I)\ |  \ Pg\in  L^2(I)\ \}
$$
and, for $g\in D^*_{min}$, $P^*_{min}g=Pg$.
\end{lemma}

\Pf
Let $g\in D^*_{min}$.  By definition, there exists $h\in L^2(I)$ such that, for all $f\in D_{min}$,
$$
\int_I Pf \bar{g} =\int_I f \bar{h}.
$$
Therefore, the distributional action of $P$ on $g$ may be identified to $Pg=h\in L^2(I)$ and we conclude that $g\in D_{max}$. 

Suppose now $g \in D_{max}= \{ \ g \in L^2(I)\ |  \ Pg\in  L^2(I)\ \}$.  By the definition of the distributional derivative, for all $f\in D'_{min}$,
 $$
 \int_{a_1}^{b_N}\bar{g}P f=  \int_{a_1}^{b_N}f\overline{Pg}
$$
so $|
 \int_{a_1}^{b_N}\bar{g}P f |\leq C|| f ||_2 
$
and we conclude that $g\in D^*_{min}$.
\EPf

Observe that  $D_{Min}\subsetneq D_{min}$.  Indeed, they have different adjoints, 
$$D^*_{min}=D_{max}\subsetneq D_{Max}=D^*_{Min},$$ by the previous proposition.

\begin{lemma} \label{lemma:py'-AC}
If $y\in D_{max}$ 
then $py'$ is absolutely continuous in $I$.
\end{lemma}
\Pf 
If $y\in D_{max}$ then $
Py=-(py')'+qy \in L^2(I)$  and therefore $(py')'\in  L^2(I)$ (and therefore in $L^1(I)$) which implies that  $py'$ is absolutely continuous.
\EPf

A more general result is the following.  In order to obtain the previous lemma, apply it to $z=1$, the constant function.

\vspace{.5cm}


\begin{prop} Consider the operator $Py=-(py')'+qy$ as above.  Then
$$
D_{max}= \{ \ y\in D_{Max}\ \vert \ py' \ \mbox{is absolutely continuous on I}\ \}.
$$
\end{prop}

\Pf
Lemma \ref{lemma:py'-AC} proves that $D_{max}\subset  \{ \ y\in D_{Max}\ \vert \ py' \ \mbox{is absolutely continuous on I}\ \}$.  On the other hand, if $y\in D_{Max}$ and $py'$ is absolutely continuous we obtain that
${Py}_{I_k}\in L^2(I_k)$ and the action of $P$ on $y$, considered as a distribution,  has a representative in $L^2(I)$ because $(py')'$ is a distribution with no atoms.

\EPf

\vspace{1cm}

\begin{rmk}
Suppose  $g \in D_{max}= \{ \ g \in L^2(I)\ |  \ Pg\in  L^2(I)\ \}$.  Then, for all continuous $f\in D_{min}$ (which is a dense subset), by the Lagrange identity
 we obtain
$$
\int_{a_1}^{b_N}(\bar{g}P f-f\overline{Pg})=\sum_i [f,g](b_i)-[f,g](a_i)=0
$$
\end{rmk}


\begin{rmk}
Observe that on an endpoint in the LP case $py'$ vanishes.  Imposing the continuity of $py'$ on the interior singular points which are in the LC case (we assume here that that both sides of a singular point are in the same case) we have
$$
 \dim  D_{Max}/D_{max}= \mbox{number of interior endpoints in the LC case}.
$$
Indeed, an integrable function which is locally absolutely continuous on the complement of isolated points
on an interval  is absolutely continuous on the interval if and only if it is continuous at these points.

If the endpoints are not in the same case we have to distinguish the number $N_c$ of LC endpoints adjacent to LC endpoints, the number  $N_p$ of  those adjacent to LP endpoints and $N_b$ the number of boundary LC endpoints (the number of LC endpoints is then $l_c=N_c+N_p+N_b$, $N_b$ is equal to 2 or 0). We obtain
$$
\dim  D_{Max}/D_{max}= \frac{N_c}{2}+N_p.
$$
\end{rmk}

\begin{rmk} 
Recall that
$$
\dim  D_{Max}/D_{Min}= 2.l_c=2(N_c+N_p+N_b).
$$
Therefore 
$$
\dim  D_{max}/D_{Min}=\dim  D_{Max}/D_{Min}-\dim  D_{Max}/D_{max}
$$
$$
=2(N_c+N_p+N_b)- (\frac{N_c}{2}+N_p)=\frac{3N_c}{2}+N_p+N_b.
$$

If there exists a self-adjoint extension $S$ of $P_{min}$ then it is also a self-adjoint extension of $P_{Min}$.
Therefore $\dim D(S)/D_{Min}=l_c$ and $\dim D_{Max}/D(S)=l_c$.  As 
$$
\dim  D_{Max}/D(S)=\dim  D_{Max}/D_{max}+\dim  D_{max}/D(S),
$$
we obtain that $\dim  D_{max}/D(S)=l_c-(\frac{N_c}{2}+N_p)=\frac{N_c}{2}+N_b$.
We conclude then that 
$$
\dim  D_{max}/D_{min}=2\dim  D_{max}/D(S)=N_c+2N_b.
$$

Note that if $N_p=0$ then  $P_{min}$ is not essentially self-adjoint only if $l_c\neq 0$
\end{rmk}

\subsection{The circle case}

The following corollary follows from the above result:

\begin{cor}\label{corollary:circle}
Let $P$ be a Sturm-Liouville operator defined on $L^2(S^1)$,
Write $S^1$ as a graph,
$$
S^1=\bigcup_{k=1}^{N}\overline{I_k},
$$
 such that all endpoints of edges  are singular points of the same nature (LP or LC) at both sides.  
Let $n_c$ be the number of  $LC$ points.
Then
$$
 def (P_{Min})= \sum_{k=1}^{N} def({P_k}_{Min})=2n_c,
 $$
 and
 $$
 def (P_{min})= n_c.
 $$
 
\end{cor}

\vspace{1cm}

{\bf{Remark}}: 
\begin{enumerate}
\item If all singular points are  $LC$ we obtain 
$$
 def (P_{Min})= 2\times \mbox{number of singular points}.
$$
\item 
If all points are  $LP$ we obtain 
$
 def (P_{Min})=0
$
so that $P_{Min}$ is self-adjoint.
\end{enumerate}

%
%
%
\subsection{Boundary conditions}

Here we state the descriptions of self-adjoint extensions using boundary conditions following 
Glazman-Krein-Naimark (see \cite{N} chapters 17-19) and the generalization for the multi-interval case (finite or infinite) by \cite{GK} and \cite{EZ1}.  

We fix again $\mathcal{I}=\bigcup_{k=1}^{N}I_{I_k}$ as before and
consider the closed operator $P$ defined on $D_{Min}$ and its adjoint $P^*$ defined on $D_{Max}$ ($N$ could be $\infty$).

\begin{dfn}
For $y,z\in D_{Max}$ define the hermitian bilinear form on $D_{Max}$ by
$$
l(y,z)=\sum \left([y,z](b_k)-[y,z](a_k)\right)
$$
\end{dfn}
This definition makes sense even for an infinite number of intervals as explained in \cite{GK} and \cite{EZ1} (the sum is absolutely convergent).

The following characterization of $D_{Min}$ paves the way to the description of boundary conditions.
\begin{thm}For any $y,z\in D_{Max}$ we have
$$
\int_I(P^*y\bar{z}-y\overline{P^*z})=l(y,z).
$$
Also,
$$
D_{Min}=\{\ y\in D_{Max}\ | \ l(y,z)=0 \ \mbox{for all}\ z\in D_{Max}\ \}.
$$
\end{thm}

%

\begin{dfn}
\begin{enumerate}
\item
Suppose $def(P_{Min})=d\in \N$.  A set of vectors $\{ \beta_i, \ 1\leq i\leq d\ \}$ in $D_{Max}$ is a generalized boundary condition set if it is linearly independent in $D_{Max}/D_{Min}$ and 
$l(\beta_i,\beta_j)=0$ for all $i,j$.   
\item If $\dim D_{Max}/D_{Min}=\infty$, the set $\{ \beta_i, \ i\in  \N\ \}$ in $D_{Max}$ is a generalized boundary condition set if it is linearly independent and if
 $\beta\in \langle \beta_i , i\in \N\rangle ^\perp$ and $l(\beta, \beta_i)=0$ for all $i\in \N$ implies $\beta\in D(P_{Min})$.
 \end{enumerate}
 \end{dfn}

 \begin{thm}\label{theorem:GKN}
 Suppose $def(P_{Min})=d$ (here $d\in \N\cup \{\infty\}$) .  Let $ \{ \beta_i \}$ be a generalized boundary condition set.
Then the operator $S$ defined by
$$
D(S)= \{\ y\in D_{Max}\ | \ l(y,\beta_i)=0 \ \mbox{for all}\  i\ \},
$$
and $Sy=P_{Max}y$
is a a self-adjoint extension of $P_{Min}$.  Conversely, if $S$ is a self-adjoint extension there exists a generalized boundary set so that $S$ is defined as above.
 \end{thm}

\section{Singular Sturm-Liouville operators on  graphs}\label{section:SLoperators}

Regular differential operators on graphs were studied in \cite{C}.  In this section we extend our results of the previous sections  to 
singular Sturm-Liouville operators on a graph.

We let $\mathcal{G}$ be a graph consisting of a set $\mathcal{V}$ of vertices and a set (finite or infinite countable) of directed edges 
in $\mathcal{V}\times \mathcal{V}$.   Each edge is identified to an interval  $[a_n,b_n]\subset \R$ and one writes
$a_n\sim v$ if the endpoint $a_n$ is identified to the vertex $v$ (by abuse of notation,  we  denote an edge by the corresponding interval).  The metric and measure on each edge induces a metric and a measure on the graph.  We have that $L^2( \mathcal{G})=\oplus_nL^2([a_n,b_n])$.

Consider $p$ and $q$ continuous functions on $\mathcal{G}$ which are smooth restricted to each interval $[a_n,b_n]$.  Moreover, we suppose that $p$ does not vanish on the interior of each edge.  
As in previous sections we define operators by choosing appropriate domains and letting $P$ act on $L^2(\mathcal{G})$ considered as a distributions space.  

\begin{dfn}
\begin{enumerate}
\item $P_{Max}$ is  the operator defined on the domain 
$$
D_{Max}=\{ \ y \in L^2( \mathcal{G})\ |  \ P{y_|}_{[a_n,b_n]}\in L^2([a_n,b_n])\ \}.
$$
\item  
$P_{Min}$ is the closure of the operator $P$ defined on the domain 
$$
D_{Min}=\{ \ y \in L^2( \mathcal{G})\ |  \ {y_|}_{[a_n,b_n]}\in C^\infty_0((a_n,b_n))\ \}.
$$

\end{enumerate}
\end{dfn}

As before we have the following result for operators on direct sums:

\begin{lemma} For the operators $P_{Min}$ and $P_{Max}$ defined as above we have
\begin{enumerate}
\item $P^*_{Min}=P_{Max}$ and $P^*_{Max}=P_{Min}$
\item $P^*_{Min}$ is  densely defined and its deficiency indices are
$$
def^+=def^-=\sum_{k=1}^N def_k,
$$
where $def_k$ is the deficiency index of each ${P_k}_{_{Min}}$, the restriction of $P$ to functions in $L^2([a_k,b_k])$.
\end{enumerate}
\end{lemma}

The proof of theorem B now follows from the formula above and the criteria in theorem \ref{theorem:deg-LCP}.  In particular,  we conclude that $P_{Min}$ is essentially self-adjoint if and only if all ${P_k}_{_{Min}}$ are essentially self-adjoint and,
therefore, if and only if all Hamiltonian trajectories are complete in $T^*(\mathcal{G}\setminus \mathcal{V})$.  This happens in the case that $p$ vanishes at each vertex with order greater than or equal to 2.  More generally, the operator $P_{Min}$ is essentially self-adjoint if and only if all operators ${P_k}_{_{Min}}$ have endpoints in the limit point case (even if the functions $p$ and $q$ have less regularity and there is no Hamiltonian field).  Observe also that if an endpoint $c$ is in the LP case than for any $y\in D_{Max}$, $py'(c)=0$.

If there exists endpoints which are regular or in the LC case there exists a space of self-adjoint extensions of the operator $P_{Min}$ which are parametrized by a unitary group (infinite if the number of endpoints in the  LC case is infinite).  Instead of considering all self-adjoint extensions, it is natural to restrict to self-adjoint extensions with certain properties. Locality of boundary conditions is a common  physical requirement. 



We define a local generalized boundary condition as a set of functions $\{ \beta_i, \ 1\leq i\leq d\ \}$ in $D_{Max}$ which is a generalized boundary condition set  (that is, it is linearly independent in $D_{Max}/D_{Min}$ and 
$l(\beta_i,\beta_j)=0$ for all $i,j$) such that each function has support around a unique vertex.
Denoting by $i_v$ the indices around the vertex $v$ and $a_{n}$ the endpoints corresponding to $v$, we have
$$
l(\beta_{i_v},\beta_{j_v})= \sum_{\ \tiny{a_n} \ \mbox{\tiny{regular or LC}}} \epsilon_n (\beta_{i_v}p\bar{\beta'}_{j_v}-\bar{\beta}_{j_v}p{\beta'_{i_v}})(a_n).
$$
Here $\epsilon_n=1$ or $-1$ depending on whether $a_n$ is the left or right endpoint of the interval.
In the following we give natural choices  of local boundary conditions.

\vspace{.5cm}

\subsection{Basis of local boundary conditions}

\subsubsection{Regular case}

Define, for any endpoint $a_n$ which is regular, $\gamma_{a_n}, \alpha_{a_n}$ smooth functions up to the endpoints with supports in half of the interval and  $\gamma_{a_n}(a_n)=\alpha'_{a_n}(a_n)=
1$,  and $\gamma'_{a_n}(a_n)=\alpha_{a_n}(a_n)=0$.
We obtain then
$$
l(\phi,\gamma_{a_n})=  
\epsilon_n 
(\phi p\bar{\gamma}_{a_n}'-
\bar{\gamma}_{a_n}p{\phi'})(a_n)=-\epsilon_n 
p{\phi'}(a_n)
$$
and
$$
l(\phi,\alpha_{a_n})=  
\epsilon_n 
(\phi p\bar{\alpha}_{a_n}'-
\bar{\alpha}_{a_n}p{\phi'})(a_n)=-\epsilon_n 
p{\phi}(a_n),
$$

\vspace{.5cm}

\subsubsection{ LC case}

Consider, for simplicity sake,  the endpoint $a_n$ at the left of the interval which is in the LC case.   In order to have canonical boundary conditions, we  choose
$u_n,v_n\in D_{Max}\setminus D_{Min}$, with support on half of the interval $[a_n,b_n]$ such that 
$$
u_n(a_n)=1,\ u_n'(a_n)=0 \ \ \mbox{and}\ \ pv_n'(a_n)=1,\ Pv_n=0
\ \mbox{near the endpoint}\ a_n.
$$

Then 
$$
l(\phi,u_n)=  
p{\phi'}(0)
$$
and
$$
l(\phi,v_n)=  
\lim_{x\to 0}\left(\phi(x)- p(x).v_n(x) .{\phi'}(x)\right)
$$

{\bf Example}: In the interval $[-L,L]$ we consider the Legendre operator 
$$
Py=((L-x^2)y')'.
$$
Let $u$ and $v$ be solutions of $Py=0$ given by :
$$
u\equiv 1,\ \  v=-\frac{1}{2L}\ln{\frac{L-x}{L+x}}
$$
Here, the left endpoint is $a=-L$ and we have
$$
l(\phi,u)=  
\lim_{x\to 0}(L-x^2){\phi'}
$$
and
$$
l(\phi,v)=  
\lim_{x\to 0}\left(\phi(x)+\frac{1}{2L}(L-x^2)\ln{\frac{L-x}{L+x}}{\phi'}(x)\right)
$$

\vspace{1cm}
%
%
%
%
%

%

\subsubsection{Singular quantum graphs with local boundary conditions}\label{section:singular-quantum-graph}

In this section we assume that the graph is finite.
A natural generalization of quantum graphs to the singular case is the following choice of 
a local boundary condition.

Vertices are partitioned into singular and regular, $\mathcal{V}_s$ and $\mathcal{V}_r$ respectively.

\begin{dfn}
\begin{enumerate}
\item
 For each vertex $v\in \mathcal{V}_r$ and $a_i\sim v$, $2\leq i\leq n$, define 
 $$
 \beta_{v_i}=\alpha_{a_i}-\alpha_{a_1}
 $$ 
and
$$
\beta_v=\sum_{ a_k\sim v} \epsilon_k \gamma_{a_k}
$$
where $\epsilon_k=1$ or $-1$ depending on whether $a_k$ is the right or left endpoint of the interval.

\item
 For each vertex $v\in \mathcal{V}_s$ and $a_i\sim v$, $2\leq i\leq n$, define 
 $$
 \beta_{v_i}=v_{a_i}-v_{a_1}
 $$ 
and
$$
\beta_v=\sum_{ a_k\sim v\  {\tiny{\mbox{is LC} }}} \epsilon_k u_{a_k}
$$
where $\epsilon_k$ is as before.
\end{enumerate}
The local boundary condition defined by  the set of functions $B=\{\beta_{v_i},\beta_{v}, \ v\in \mathcal{V}\}$ is called  (singular) Neumann-Kirchhoff  boundary  condition.
\end{dfn}

 One can state then a particular self-adjoint extension:

\begin{thm}\label{theorem:singular-quantum-graph}
The extension of $P_{Min}$ defined by  $P_{Max}$ restricted to the domain 
 $$
\{ \ y \in D_{Max}\ | \ B(y)=0 \ \mbox{for all}\ v\in \mathcal{V}\ \}
$$
 is self-adjoint.
\end{thm}

\Pf 
In the regular vertex case, as
$
l(\alpha_{a_m},\alpha_{a_n})=0,
$
we obtain that  $l(\beta_{v_i},\beta_{w_j})=0$ for all $v,w\in \mathcal{V}$.  Also
$l(\gamma_{a_m},\alpha_{a_n})=0$ which implies that
$$
l( \beta_{v_i},\beta_w)=0,
$$
for all $v,w$ and $i$.  The case of a singular vertex is similar.

One chooses the $N=def(P_{Min})$ functions $ \beta_{v_i},\beta_v$ in $D_{Max}$.  They are clearly linearly independent and satisfy the condition   By GKN extension theorem (Theorem \ref{theorem:GKN}),  the set of functions
$$
D= \{\ y\in D_{Max}\ | \ l(y,\beta_{v_i})=0, l(y,\beta_{v})=0 \ \mbox{for all}\  v\ \}
$$
is the domain of a self-adjoint extension of $P_{Min}$.

\EPf

\begin{rmk}
In the case that all vertices are regular, the conditions $l(y,\beta_{v_i})=0$ imply that the 
the function $y$ is continuous on $\mathcal{G}$ so
$$
D(S)= \{\ y\in D_{Max}\ | \ y \ \mbox{is continuous and}\ l(y,\beta_{v})=0, \mbox{for all}\  v\ \}.
$$

In particular, if $\beta_{v}=\sum_{ a_k\sim v} \epsilon_k \beta_{a_k}$ (where $\epsilon_k\in \{-1,1\}$ as before) then the above boundary condition coincides with Neumann-Kirchhoff boundary condition of quantum graphs. This is the natural boundary condition which makes regular vertices of valence two disappear.

\end{rmk}

%

%

\subsubsection{Remark: Local self-adjoint extensions on a graph}\label{reamrk:delta}

An analog to the $\delta$-boundary conditions in quantum graphs is obtained from the basis described above.  More explicitly for LC endpoints we impose the following conditions.  For a fixed  $\alpha\in \R$,

\begin{enumerate}
\item For every edge $e$ abutted to a vertex $v\in {\mathcal{V}}$ satisfying $Pv=0$ and $pv=0$ at the endpoints, the function
$$
\lim_{t\to v}(y- v py')
$$
has the same limit.  That is, the function $y- v py'$ is a continuous function on the graph.
\item For every vertex $v\in {\mathcal{V}}$ and $a_i\sim v$,
$$
\sum_{ a_k\sim v\  {\tiny{\mbox{is LC} }}} \epsilon_k py'(a_k)=\alpha (y- v py')(v)
$$
where $\epsilon_k$ is as before.
\end{enumerate}


\subsection{The spectrum}\label{section:spectrum}

\subsubsection{Interlacing}
\label{subsection:interlacing}

In order to obtain Weyl's law for a singular quantum graph we will need an interlacing estimate between self-adjoint extensions with the canonical Friedrichs 
extension.  We will quote the following result, which is a special case of theorem 1.2 in  \cite{BGLS} et all.

\begin{thm} \label{theorem:interlacing}
Let $P_m$ be the minimal operator associated to a singular quantum graph (which is bounded from below), $P_F$ be the Friedrichs extension and $P$ any self-adjoint extension.  Then, for any $[a,\lambda]\in \R\setminus Spec_{ess}P_m$,
$$
|N(P,(a,\lambda])-N(P_F,(a,\lambda])|\leq def(P_m).
$$
\end{thm}

In particular if all the singular points are LC and  the operator is given by $P_my=-(py')'+qy$ with $p>0$ in the interior of each edge, we obtain 
$$
|N(P,(-\infty,\lambda])-N(P_F,(-\infty,\lambda])|\leq def(P_m).
$$

\subsubsection{Spectrum of Friedrichs extensions on an interval}\label{subsection:friedrichs-extension}

We consider the simplest operators with LC or regular endpoints.   There are three classical cases:

\begin{enumerate}
\item The free Schr\"odinger operator on $(0,L)$:
$$
P_m y=-y'',
$$
with both endpoints regular.  The Friedrichs extension is the Dirichlet operator with spectrum
$$
\lambda_n=\left(\frac{n\pi}{L}\right)^2
$$
for $n\geq 1$.  We obtain
$$
\sqrt{\frac{\lambda L}{\pi}}-1 \leq N(\lambda,P_F)\leq \sqrt{\frac{\lambda L}{\pi}}.
$$
\item The Legendre operator on $(-L/2,L/2)$:
$$
P_m y=-(\frac{L^2}{4}-t^2)y'',
$$
with both endpoints LC (non-regular). The Friedrichs extension is the operator with boundary conditions $\lim_{t\to \pm L/2}{(1-t^2)y'}=0$. The spectrum is 
$$
\lambda_n=n(n+1),
$$
for $n\geq 0$, independent of $L$ with eigenfunctions, the Legendre polynomials $P_n(2t/L)$.
Note that it does not depend on the size of the interval.  The counting function may be estimated:
$$
\sqrt{\lambda}-2 \leq N(\lambda,P_F)\leq \sqrt{\lambda}+2.
$$
\item The equation on $(0,L)$:
$$
P_m y=-(ty')',
$$
with $0$ an LC point and $L$ a regular point. Eigenvalues are obtained by solving $-(ty')'=\lambda y$.   By a change of variable $z = 2\sqrt{\lambda t}$ one obtains
the Bessel equation of order 0: $z^2 Y'' + z Y' + z^2 Y = 0$ where $Y(z)=y(t)$.  Therefore
$$
y(t)= AJ_0(2\sqrt{\lambda t}) + BY_0(2\sqrt{\lambda t}),
$$
where $J_0$ and $Y_0$ are Bessel functions.
The Friedrichs extension is given by the boundary conditions $y(L)=0$ and $\lim_{t\to 0}{ty'}=0$.
The eigenvalues of the Friedrichs extension are given by 
$$
\lambda_n=\frac{j_{0,n}^2}{4L}
$$
for $n\geq 1$.  Here $j_{0,n}$ are the roots of the Bessel function $J_0$.  For large $n$, we have
$$
\lambda_n\sim \frac{n^2\pi^2}{4L}.
$$
We obtain, for large $\lambda$, 
$$
 N(\lambda,P_F)\sim \sqrt{\frac{2}{\pi}}L^{1/4}\sqrt{\lambda}.
$$

\end{enumerate}

\subsubsection{Weyl's law}

In this last section we establish a Weyl's formula for a Sturm-Liouville graphs which generalizes the formula in the case of regular quantum graphs (see \cite{BK}).

\begin{thm}\label{theorem:Weyl}
Let $(\mathcal{G},P_{min})$ be a Sturm-Liouville graph with vertices which are either regular or  LC. We suppose that the principal symbol is positive in each interval. Then for any self adjoint extension $P$ of $P_{min}$, we have
$$
 N(\lambda,P)\sim \left(\frac{1}{\sqrt{\pi}}\left( \sum_{e\in R}{L_e}^{1/2}
 +{\sqrt{2}}\sum_{e\in RLC}{L_e}^{1/4}\right) +
  \#\{\ e\in E\ |\ e\in E_{LC}\ \}
  \right)\sqrt{\lambda},
$$
where the first sum is made over the edges with both endpoints regular, the second with one regular  endpoint and the other LC and the last term is the number of edges with both endpoints LC (non-regular).
\end{thm}

\Pf
By the interlacing theorem we have that if there are $N$ edges in the graph
$$
| N(\lambda,P)- N(\lambda,P_F)|\leq N,
$$
where $P_F$ is the Friedrichs extension.
It suffices therefore to estimate $N(\lambda,P_F)$. 

The Friedrichs extension disconnects the graph.  Its spectrum is the spectrum of a union of edges and the counting functions of edges add up to obtain the total counting function.
We have
$$
N(\lambda,P_F)\sim C\sqrt{\lambda}
$$
with the constant
$$
C=\frac{1}{\sqrt{\pi}}\left( \sum_{e\in R}{L_e}^{1/2}+{\sqrt{2}}\sum_{e\in RLC}{L_e}^{1/4}\right) + \#\{\ e\in E\ |\ e\in LC\ \}.
$$
\EPf

{\bf{Example:  Legendre graphs}}:
Consider a finite graph where in each edge $e_i$ edges we impose the Legendre operator 
$$
{P_|}_{e_i}y=-((L_i^2-x^2)y')'.
$$
Weyl's law becomes, for every self-adjoint extension, 
$$
 N(\lambda,P)\sim 
  \#\{\ e\in E\ |\ e\in LC \ \}\sqrt{\lambda}.
$$
Note that, for Legendre's graphs, the constant in Weyl's law does not depend on the metric.

\section{Appendix}

\subsection{Criteria for endpoints in the limit point or limit circle case}

In this appendix we prove Theorem \ref{theorem:deg-LCP}.
First we prove a related proposition. It uses  Weidemann Theorem 6.3 which, for the sake of completeness, we 
detail a proof.

%
%
%

\begin{prop} 
The  endpoint $b$  (respectively $a$) is a limit point case  if for  $c\in I$, $\frac{b-t}{\sqrt{p(t)}}\notin L^1(c,b)$ (respectively $\frac{t-a}{\sqrt{p(t)}}\notin L^1(a,c)$).
\end{prop}

\Pf
\begin{enumerate}
\item
First we prove the following (see Weidemann Theorem 6.3):  $b$ is a limit point case if, for some $c\in I$,
$$
\int _c^x\frac{1}{{p(t)}}dt\notin L^2(c,b).
$$

Choose $x_0\in I$ such that, for all $x\in [x_0,b)$,  $q(x)>\gamma$.
Then the solutions to the  equation $(P-\gamma)u=0$ with initial conditions $u(x_0)=1, u'(x_0)=0$ or $u(x_0)=0, p(x_0)u'(x_0)=1$ satisfy
$$
(p(x)u'(x))'= (q(x)-\gamma)u(x).
$$
Therefore $p(x)u'(x)$ and $u(x)$ are increasing in $[x_0,b)$.  Choose now $c\in (x_0,b)$ so that,
for $x\in (c,b)$, $p(x)u'(x)\geq p(c)u'(c)>0$.

We obtain then 
$$
u(x)=u(x_0)+\int_{x_0}^xu'(t)dt>u(x_0)+p(c)u'(c)\int_{x_0}^x\frac{ 1}{p(t)}dt.
$$
Therefore, if $\int _c^x\frac{1}{{p(t)}}dt\notin L^2(c,b)$ then $u\notin L^2(c,b)$.  We conclude that the endpoint $b$ is in the limit point case.

\item
We show now that $\frac{1}{\sqrt{p}}\notin L^1(c,b)$ implies $\int _c^x\frac{1}{{p(t)}}dt\notin L^2(c,b)$.
Indeed, from Cauchy-Schwartz inequality, we have
$$
\left ( \int _c^x\frac{1}{\sqrt{p(t)}}dt \right )^2\leq
\int _c^x\frac{1}{{p(t)}}dt \int_c^x dt\leq (b-c)\int _c^x\frac{1}{{p(t)}}dt
$$
So,
$$
\int_c^b \left ( \int _c^x\frac{1}{\sqrt{p(t)}}dt \right )^4 dx\leq
(b-c)^2\int_c^b\left (\int _c^x\frac{1}{{p(t)}}dt\right )^2 dx.
$$
Observe now that, for $x$ sufficiently close to $b$, as $\frac{1}{\sqrt{p}}\notin L^1(c,b)$ so
 $$
  \int _c^x\frac{1}{\sqrt{p(t)}}\,dt \leq  \left ( \int _c^x\frac{1}{\sqrt{p(t)}}dt \right )^4
 $$
  and therefore
$$
\int_c^b \int _c^x\frac{1}{\sqrt{p(t)}}\,dt  \, dx
\leq
(b-c)^2\int_c^b\left (\int _c^x\frac{1}{{p(t)}}dt\right )^2 dx.
$$
This implies that 
$$
\int_c^b \int _c^b\theta(x-t)\frac{1}{\sqrt{p(t)}}\,dt  \, dx
\leq
(b-c)^2\int_c^b\left (\int _c^x\frac{1}{{p(t)}}dt\right )^2 dx,
$$
where $\theta(x-t)=0$ if $x<t$ and otherwise it is identically equal to $1$.
We obtain, using Fubini,
$$
\int_c^b \frac{b-t}{\sqrt{p(t)}}\,dt  
\leq
(b-c)^2\int_c^b\left (\int _c^x\frac{1}{{p(t)}}dt\right )^2 dx.
$$
This inequality proves the proposition.

\end{enumerate}
\EPf

\subsubsection{Proof of Theorem \ref{theorem:deg-LCP}}

The theorem may be compared to Theorem 6.4 in Weidemann.  It is probably a classical result but we give a proof for the sake of completeness. The easiest implication is given in the first item.  It gives a more useful criterium than the previous proposition for an endpoint to be in the limit point case 
when the principal symbol is of the form $p(x)=(b-x)^n\phi(x)$ where $\phi$ is $C^1$ up to the boundary.

\vspace{.5cm}

\begin{enumerate}
\item
We evaluate for $c\in I$ sufficiently close to $b$ and $n\geq 2$: 
$$
\int _c^x\frac{1}{{p(t)}}dt=
\int _c^x\frac{1}{(b-t)^n\phi(t)}dt
 \geq 
K\int _c^x\frac{1}{(b-t)^n}dt=
K(\frac{1}{(b-x)^{n-1}}-\frac{1}{(b-c)^{n-1}}),
$$
where $K^{-1}=\min \{\  \phi(x)\ | x\in [c,b]\ \}$. Therefore, we verify that $\int _c^x\frac{1}{{p(t)}}dt \notin L^2(c,b)$ and we may apply the first item of the proof of the above proposition and conclude that the endpoint $b$ is in the limit point case.

\item 

Suppose now  $p(x)=(b-x)\phi(x)$ where $\phi$ is a continuous function up to $b$ and $\phi(b)> 0$.   We want to show that $b$ is in the limit circle case.  That is, all solutions to
$Pu=0$ are in $L^2(I)$.  We adapt an argument  of Weidemann, part b of Theorem 6.4.

We give the argument at the endpoint $a$ which we suppose $a=0$. That is, near $0$, we have $p(x)=x\phi(x)$. This simplifies notation.  The equation is then
$$
-x\phi(x)u''(x)-(\phi(x)+x\phi'(x))u'+q(x)u=0.
$$
We write it as
$$
-xu''(x)-(1+\frac{x\phi'(x)}{\phi(x)})u'+\frac{q(x)}{\phi(x)}u=0
$$
and note that $0$ is a regular singular point.

Let $u_0$ be a solution of the equation defined on $(a,c)$
$$
{P_0}u=-t u''-Lu'+Ku=0,
$$
 where $K\geq 0$ is a constant such that $K\geq \frac{q(t)}{\phi(t)}$ and $L\geq 0$ is a constant such that $1+\epsilon\geq L\geq 1+\frac{q(t)}{\phi(t)}$, for all $a\leq t \leq  c$ and $\epsilon$ sufficiently small.  Clearly,  the conditions are satisfied if $c$ is sufficiently close to $a$.

 \vspace{.5cm}
 {
 {This is an analytic differential equation with a regular singular point at $0$. It is well known that all solutions are in $ L^2(a,c)$ where $c\in (a,b)$.   Indeed, the  indicial equation is
 $$
 r(r-1)+Lr=0
 $$
 and therefore $r=0$ or $r=1-L\geq -\epsilon$.  The solutions to the equation are seen to be in $L^2(a,c)$.

   The point now is to compare solutions to the original equation to the one above.
 }
 }
  \vspace{.5cm}
 
 Consider now a solution $u$ of the equation $Pu=0$ on $(a,c)$ and let $u_0$ a solution of $P_{u_0}$ with initial conditions
 $$
 u_0(c)=| u(c)|+1,\ \  {u_0}'(c)=-|u'(c)|-1.
 $$
 
  \vspace{.5cm}
 Claim: $|u(x)|< u_0(x)$ for $x\in (a,c)$.
  \vspace{.5cm}
  
 Define
 $$
 d=\inf \{\ s\in (a,c)\ | \ u_0(x)> |u(x)|\ \mbox{and}\ u'_0(x)<- |u'(x)|\ \mbox{for}\ x\in (s,c)\ \}.
 $$ 
Suppose, by contradiction, that $d>0$.  Then
$$
u_0(d)=u(d)\ \mbox{or}\ u'_0(d)=- |u'(d)|.
$$
But 
$$
u'_0(d)= u'_0(c)+\int_c^d u_0''(x)dx
=-|u'(c)|-1-\int_d^c u_0''(x)dx
$$
$$
=-|u'(c)|-1-\int_d^c \frac{1}{x}(-Lu_0'(x)+Ku_0(x))dx
$$
$$
<  -|u'(c)|-1-\int_d^c \frac{1}{x}(|Lu'(x)|+|Ku|)dx
$$
$$
 \leq
 -|u'(c)|-1-\int_d^c \frac{1}{x}\left(|(1+\frac{p'(x)}{\phi(x)})u'(x)|+|\frac{q(x)}{\phi(x)}u(x)|\right)dx
$$
$$
\leq  -|u'(c)|-1-\int_d^c \frac{1}{x}\left(|(1+\frac{p'(x)}{\phi(x)})u'(x)+\frac{q(x)}{\phi(x)}u(x)|\right)dx=-|u'(c)|-1-\int_d^c |u''(x)|dx
$$
$$
\leq  -|u'(c)|-1-\left |\int_d^c u''(x)dx\right |= -|u'(c)|-1-\left | u'(c)-u'(d)\right |< -|u'(d)|.
$$
Also, 
$$
u_0(d)=u_0(c)+\int_c^d u_0'(x)dx= |u(c)|+1-\int_d^c u_0'(x)dx
$$
$$
>|u(c)|+1+\int_d^c |u'(x)|dx\geq |u(c)|+1+\left | \int_d^c u'(x)dx\right |= |u(c)|+1+\left | u(d)-u(c)\right |>u(d).
$$
The inequalities $u'_0(d)< -|u'(d)|$ and $u_0(d)>u(d)$ provide the contradiction.

\end{enumerate}
\EPf

\bibliographystyle{alpha}
\bibliography{biblio-operators-circle}

\begin{flushleft}
  \textsc{E. Falbel\\
  Institut de Math\'ematiques \\
  de Jussieu-Paris Rive Gauche \\
CNRS UMR 7586 and INRIA EPI-OURAGAN \\
 Sorbonne Universit\'e, Facult\'e des Sciences \\
4, place Jussieu 75252 Paris Cedex 05, France \\}
 \verb|elisha.falbel@imj-prg.fr|
 \end{flushleft}
 
\end{document}